\documentclass[12pt]{article}
\usepackage[utf8]{inputenc}
\usepackage[T1]{fontenc}
\usepackage{lmodern}
\usepackage{geometry}
\usepackage{amsmath,amssymb,amsfonts,amsthm,mathtools}
\usepackage{hyperref}
\usepackage{enumitem}
\setlist{nosep}
\usepackage{microtype}
\usepackage{titlesec}
\usepackage{graphicx}
\usepackage{float} 
\titleformat{\section}[hang]{\normalfont\Large\bfseries}{\thesection\quad}{0pt}{}
\titleformat{\subsection}[hang]{\normalfont\large\bfseries}{\thesubsection\quad}{0pt}{}

\newtheorem{theorem}{Theorem}[section]
\newtheorem{lemma}{Lemma}[section]
\newtheorem{conjecture}{Conjecture}[section]
\newtheorem{definition}{Definition}[section]

\usepackage{bm}

\title{Detecting Causality with Conjugation Quandles over Dihedral Groups}
\author{Zining Fan\\[0.5em]
Emmaus High School\\
\texttt{Zncubefan@gmail.com}}
\date{August 21, 2025}

\begin{document}
\thispagestyle{plain}
\begingroup
\centering
{\LARGE\bfseries Detecting Causality with Conjugation Quandles over Dihedral Groups\par}
\vspace{0.75em}
{\large Zining Fan\par}
\vspace{0.25em}
{\normalsize Emmaus High School, USA\par}
\texttt{Zncubefan@gmail.com}\par
\vspace{0.75em}
{\normalsize September 1, 2025\par}
\par
\endgroup

\bigskip
\begin{center}\textbf{Abstract}

\begin{minipage}{0.92\textwidth}
\bigskip
\noindent We study whether quandle colorings can detect causality of events for links realized as skies in a $(2+1)$-dimensional globally hyperbolic spacetime $X$. Building off the Allen--Swenberg paper in which their $2$-sky link was conjectured to be causally related, they showed that the Alexander--Conway polynomial does not distinguish that link from the connected sum of two Hopf links, corresponding to two causally unrelated events. We ask whether the Alexander--Conway polynomial together with different types of quandle invariants suffice. We show that the conjugation quandle of the dihedral group $D_5$, together with the Alexander--Conway polynomial, does distinguish the two links and hence likely does detect causality in $X$. The $2$-sky link shares the same Alexander--Conway polynomial but has different $D_5$ conjugation-quandle counting invariants. Moreover, for $D_5$ the counting invariant alone already separates the pair, whereas for other small dihedral groups $D_3$, $D_4$, $D_6$, $D_7$ even the enhanced counting polynomial fails to detect causality. In fact we prove more that the conjugation quandle over D5 distinguishes all the infinitely many Allen-Swenberg links from the connected sum of two Hopf links. These results present an interesting reality where only the conjugation quandle over $D_5$ coupled with the Alexander--Conway polynomial can detect causality in $(2+1)$ dimensions. This results in a simple, computable quandle that can determine causality via the counting invariant alone, rather than reaching for more complicated counting polynomials and cocycles.

\end{minipage}

\bigskip
\noindent\textbf\sloppy{Key Words:} Knot Invariant, Globally Hyperbolic Spacetime, Dihedral Group, Quandle invariants, Alexander--Conway polynomial
\medskip

\end{center}
\bigskip
\bigskip

\section{Introduction}

\subsection*{Knot and links}
A \textit{knot} is an embedding of $S^1$ into $3$-dimensional Euclidean space $\mathbb{R}^3$. These knots are considered equivalent when a sequence of Reidemeister moves can deform one knot into the other. The simplest knot is the \textit{unknot}, also called trivial knot, which is a circle in $\mathbb{R}^3$. When one or more knots are placed together in $\mathbb{R}^3$, we get a \textit{link}. The simplest link is the \textit{trivial link}, which is the union of $2$ unlinked circles in $\mathbb{R}^3$. The simplest nontrivial link is the \textit{Hopf link}, which consists of $2$ circles which are linked together.

\begin{figure}[h]
  \centering
  \includegraphics[width=0.5\textwidth]{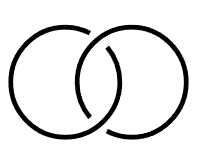}
  \caption{Hopf Link}
  \label{fig:diagram}
\end{figure}

\subsection*{Knot invariants}
A \textit{knot invariant} is defined as a function that is assigned to a knot that does not change under \textit{Reidemeister} moves. If two knots or links have different invariants, then they are not ambient isotopic, however, the converse is not true.

Knot invariants matter when studying causality in spacetimes due to the Low \cite{Low1988}, and Legendrian Low Conjecture of Nat\'ario and Tod \cite{NatarioTod2004} results of Chernov and Nemirovski \cite{ChernovNemirovski2016}: in $(2+1)$-dimensional globally hyperbolic spacetimes, causal relations between events can be detected by the linking of their skies. This allows us to translate events into link relations, and work with the relations using invariants to distinguish between link types.

When two disjoint simple longitude closed curves are embedded onto a standard torus ($T$ in $\mathbb{R}^3$) around the axial circle, the closed curves are not linked, so the corresponding skies of the two events are not topologically linked. As is shown by Low \cite{Low1988} these two curves represent skies of two causally unrelated events.

A torus can be represented as $\mathbb{R}^3$ with the trivial circle passing through the hole in a torus. We find that the transformed link in $\mathbb{R}^3$ is ambient isotopic to the connected sum of the Hopf links, $H\#H$. In extension, the connected sum of the $2$ Hopf links in $\mathbb{R}^3$ is a model for non-causally related events in a $(2+1)$-dimensional globally hyperbolic spacetime.

\medskip
\begin{figure}[h]
  \centering
  \includegraphics[width=0.6\textwidth]{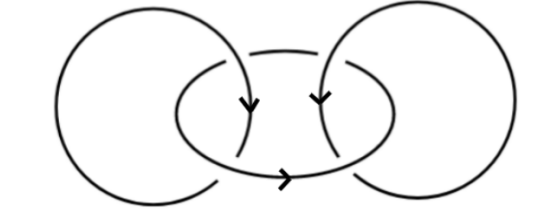}
  \caption{Connected Sum of Hopf Links}
  \label{fig:diagram}
\end{figure}

\medskip
This is where knot and link invariants become important for this research. If an invariant shows that a given link is not isotopic to $H\#H$, then that link potentially represents a causal relationship in a $(2+1)$-dimensional globally hyperbolic spacetime, since in $\mathbb{R}^3$, a link is non-causal if and only if it is $H\#H$. Of course, the link component should consist of two knots individually isotopic to longitudes in a torus since other knots are not isotopic to skies.

Nat\'ario and Tod \cite{NatarioTod2004} used a collection of links to map several causally related events and discovered that the Jones Polynomial distinguishes them from $H\#H$ thus creating Allen--Swenberg conjecture that Jones Polynomial does always capture causality.

Chernov, Martin and Petkova \cite{ChernovMartinPetkova2020} showed that more powerful invariants like the Heegaard--Floer Homology (whose graded Euler characteristic is the Alexander--Conway Polynomial) and Khovanov Homology (whose graded Euler characteristic is the Jones Polynomial) can detect causality in a $(2+1)$-dimensional globally hyperbolic spacetime with the Cauchy surface $\Sigma \neq S^2$ and $\Sigma \neq \mathbb{R}P^2$. Allen and Swenberg used empirical evidence to suggest that the Jones polynomial could detect causality in a $(2+1)$-dimensional globally hyperbolic spacetime while the Alexander--Conway polynomial could not \cite{AllenSwenberg2021}. Since for their $2$ sky link, The Jones polynomial was distinct from the connected sum of the Hopf links, while the Alexander--Conway polynomial was not.

A \textit{quandle} is a special algebraic structure that remains unchanged under Reidemeister moves. Thus it is an invariant of knots and links. Quandles were first separately introduced in 1982 by Joyce and Matveev \cite{Joyce1982,Matveev1984}. Joyce used the term quandles to describe these algebraic structures, while Matveev described them as ``distributive groupoids'', hinting at the nature of the self-distributivity of a quandle, and the bilinear group-like operation of these structures.

Recently, new research has been done on using quandles and the Alexander--Conway polynomial to detect causality. In a 2022 paper, Jack Leventhal \cite{Leventhal2023} showed that when even using together Affine Alexander quandles and the Alexander polynomial, it could not detect causality in a $(2+1)$-dimensional globally hyperbolic spacetimes. He also tried, in his paper, to see if Takasaki quandles, the $t = -1$ specialization of Alexander quandle, could detect causality. He combined this with the Alexander polynomial but still found that it failed to do so in a $(2+1)$-dimensional globally hyperbolic spacetime. This result was later augmented with a paper by Hongxu Chen \cite{Chen2024}, who found a computational mistake, but eventually came to the same conclusion.

In 2024, Ayush Jain \cite{Jain2024} followed the same approach with Symplectic Quandles over the spaces $\mathbb{Z}_3$, $\mathbb{Z}_4$, and $\mathbb{Z}_5$, and the matrix $\begin{bmatrix} 1 & 0 \\[2pt] 0 & -1 \end{bmatrix}$. He found that the number of homomorphisms i.e.\ the counting invariant could not detect causality in a $(2+1)$-dimensional globally hyperbolic spacetime, however, the \textbf{enhanced quandle counting polynomial could}. In this paper, we investigate if conjugation quandles, along with the Alexander polynomial can detect causality in a $(2+1)$-dimensional globally hyperbolic spacetime, or more precisely if the conjugation quandle does distinguish the first Allen--Swenberg link from $H\#H$, and if so, can it distinguish between all Allen--Swenberg links and $H\#H$.

\section{Globally Hyperbolic Spacetimes}

\begin{definition}[Spacetime and Lorentz dot product]\label{def:spacetime}
A \textit{spacetime} $X$ is defined as a Lorentz manifold with an operation called the Lorentz dot product
\begin{equation*}
(x_1, x_2, \dots, x_n, t_1)\cdot (y_1, y_2, \dots, y_n, t_2) = x_1 y_1 + x_2 y_2 + \dots + x_n y_n - t_1 t_2 \, .
\end{equation*}
\end{definition}
A \textit{causal trajectory} between the points (events) $p$ and $q$ from a curve $y'$, must satisfy the inequality
\begin{equation*}
y'(x,t)\cdot y'(x,t) \le 0 \, .
\end{equation*}
If so, the points are \textit{causally related}, and the Lorentz dot product of the velocity for the causal trajectory with itself must be timelike or null (slower or equal to the speed of light).

\begin{definition} Globally Hyperbolic Spacetimes \cite{BernalSanchez2007} can be defined as having a Cauchy surface $\Sigma$ which is the set of all points $p$ such that every past inextensible timelike curve through $p$ passes though it at exactly one point. While the $\mathbb{R}$ coordinate is one of many timelike functions.
\end{definition}
\begin{theorem}[Bernal and S\'anchez \cite{BernalSanchez2003}]\label{thm:bernal-sanchez}
A globally hyperbolic spacetime satisfies:
\begin{enumerate}
\item For all $p,q\in X$, $J^+(p)\cap J^-(q)$ is compact.
\item No time travel.
\end{enumerate}
\end{theorem}

The first condition is known as the \textit{absence of naked singularities}. Here for $p$ and $q$ in a differentiable manifold $M$, $J^+(p)$ is defined as the \textit{causal future} of a point $p$, which is the set of all points in the spacetime that can be reached from $p$ by a future-directed causal trajectory. On the other hand, $J^-(q)$ is defined as the \textit{causal past}, which is the set of all points which can influence $q$ via a future-directed causal trajectory. The causal diamond, which is the intersection set of all points in the causal future of $p$, and the set of all points in the causal past of $q$, is compact in a globally hyperbolic spacetime. Penrose \cite{Penrose1979}, under his strong cosmic censorship conjecture, conjectured that all physically reasonable spacetimes follow these axioms.

The space of light rays in a globally hyperbolic $X$ is homeomorphic to the spherical cotangent bundle $S(T^*\Sigma)$ of a Cauchy surface $\Sigma$. For the case $\Sigma = \mathbb{R}^2$ this is a doughnut (torus). Using standard covering techniques, as in Chernov--Martin--Petkova, we obtain that the links in the doughnut work for all globally hyperbolic $(2+1)$--dimensional spacetimes with $\Sigma \ne S^2$ and $\Sigma \ne \mathbb{R}P^2$. We then get a sphere of light rays passing through a point $x$, called the \emph{sky} of $x$. (This sphere lives in the space of all light rays rather than in a spacetime.)

\begin{conjecture}[Low Conjecture \cite{Low1988}]\label{conj:low}
Assume that the universal cover of a smooth spacelike Cauchy surface of a globally hyperbolic $(2+1)$--dimensional spacetime $(X,g)$ is $\mathbb{R}^2$. Then the skies of causally related points in $X$ are topologically linked.
\end{conjecture}

The resulting link from their skies is named to be non-trivial while the link of two longitudes corresponding to two causally unrelated events is called trivial. This conjecture was proven by Chernov and Nemirovski \cite{ChernovNemirovski2016}. Here, they showed that the relation holds as long as the Cauchy surface $\Sigma \not\simeq S^2, \mathbb{R}P^2$.

Low conjecture is only applicable to $(2+1)$-dimensional globally hyperbolic spacetimes. This is not useful for our universe's spacetime, which is $3+1$, Therefore, Nat\'ario and Tod \cite{NatarioTod2004} gave name to the Legendrian Low conjecture, which states that if the skies of $2$ events are Legendrian linked, they are causally related. This was also proven by Chernov and Nemirovski \cite{ChernovNemirovski2016}. For the generalizations of these results see also \cite{ChernovNemirovski2016}, \cite{ChernovNemirovski-Low}, \cite{ChernovNemirovski-NonNeg}, and related works.

\section{Quandles}

\subsection{Quandle properties}

\begin{definition}[Quandle]\label{def:quandle}
A quandle is defined as a set $X$ with a binary operation that satisfies the following three axioms \cite{NelsonNotes2004}:
\begin{enumerate}
\item $x \triangleright x = x$, for all $x \in X$ \textbf{(idempotency)};
\item For elements $x,y \in X$, there exists some element $z$ such that $x = y \triangleright z$ \textbf{(right invertibility)};
\item $(x \triangleright y) \triangleright z = (x \triangleright z) \triangleright (y \triangleright z)$, for all elements $x,y,z \in X$ \textbf{(self-distributivity)}.
\end{enumerate}
\end{definition}

The first property indicates the idempotency property of the quandle. The second property is the existence of the right inverse, $\triangleright^{-1}$, where $(x \triangleright y) \triangleright^{-1} y = x$. The third property is the self-distributivity of the quandle. Quandles are generally neither commutative nor distributive.

Joyce introduced the fundamental quandle knot invariant. Under the fundamental quandle, he showed this relationship between the arcs and crossings \cite{Joyce1982}.

\medskip
\begin{figure}[t]
  \centering
  \includegraphics[width=0.6\textwidth]{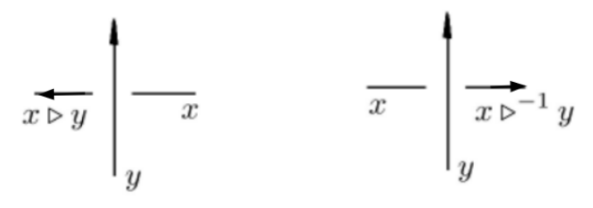}
  \caption{Crossing Relations [10]}
  \label{fig:diagram}
\end{figure}

\bigskip
\bigskip
\bigskip
\bigskip
\bigskip
\bigskip
\bigskip
Figure 3 illustrates the positive crossing, where the crossing is right-handed. Here, the operation takes the arc $x$, which passes underneath the arc $y$ from the left and combines them into $x \triangleright^{-1} y$. Diagram B shows the positive crossing, where the crossing is left-handed. The operation here takes the arc $x$, which passes underneath the arc $y$ from the right and combines it into $x \triangleright y$.

Joyce \cite{Joyce1982}, and Matveev \cite{Matveev1984}, showed that the fundamental quandle did not change over Reidemeister moves, and so it was classified as a knot (link) invariant. They also discovered that the fundamental quandle, $Q(L)$, of an oriented knot (link) is a complete knot (link) invariant and its operations are shown by the figure above. However, the fundamental quandle is very tricky to compute and mostly impossible to use because this requires showing that two quandles are not isomorphic.

\subsection{Quandle invariants}

When a quandle is applied to a link, you get the coloring for each crossing. These colorings have systems of equations that can be solved over a group. Quandle homomorphisms are simply the number of solutions to these systems of equations. These homomorphisms can be characterized by the property $f(x \triangleright y) = f(x) \triangleright f(y)$.

In a link, the quandle coloring is defined as a way to label the link with the operations above.

\begin{definition}[Counting invariant]\label{def:counting}
The quandle counting invariant, denoted by
\begin{equation*}
\phi_T(L) = \left|\mathrm{Hom}(Q(L), T)\right| ,
\end{equation*}
counts the number of homomorphisms for the link (ways that a link can be colored).
\end{definition}

If the number of homomorphisms is different for two links, then the links are different. However, if the number of quandle homomorphisms are equal, it is inconclusive whether the link is the same.

The quandle counting invariant is often insufficient for the purpose of determining causality, since most of the time, it cannot detect the difference between different knots and links. Thus, Nelson introduced the idea of an enhanced quandle counting invariant.

\begin{definition}[Enhanced counting polynomial \cite{NelsonNotes2004}]\label{def:enhanced}
The enhanced quandle counting polynomial is defined as follows:
\begin{equation*}
\Phi_E(L, T) \;=\; \sum_{f \in \mathrm{Hom}(Q(L),T)} q^{|\mathrm{Im}(f)|}\,.
\end{equation*}
We examine all quandle homomorphisms and, for each one, record how many colors appear in its image; this number is used as a weight. Collecting these counts into coefficients refines the ordinary counting invariant and can distinguish links that the unweighted count cannot.
\end{definition}

There are many useful types of quandles. The fundamental quandle is the structure determined by the relations $x \ast y$ arising at the crossings. Other examples include the affine Alexander quandle, which Leventhal \cite{Leventhal2023} showed does not detect causality, and the symplectic quandle, which Jain \cite{Jain2024} showed to detect causality.

\subsection{Conjugation quandles}

In this paper, we will be exploring conjugation quandles, and whether or not they can detect causality in our $(2+1)$ dimensional spacetime.

\begin{definition}[Conjugation quandle]\label{def:conj-quandle}
Given a group $G$, the conjugation quandle on $G$ is the set $G$ with operation
\begin{equation*}
a \triangleright b \;=\; b^{-1} a b \, .
\end{equation*}
This operation satisfies idempotency, right-invertibility, and self-distributivity.
\end{definition}

In a paper by Allen and Swenberg \cite{AllenSwenberg2021}, they showed that there were a family of links that were conjectured to be causally related, but for which causality could not be detected by the Alexander--Conway polynomial. In this paper our goal is to see whether the conjugation quandles over various dihedral groups could possibly differentiate between the connected sum of the Hopf links and the Allen--Swenberg links, so that when combined with the Alexander--Conway polynomial, be able to detect causality for some and then conjecturally for all links.

\section{Detecting Causality}

\subsection{Dihedral group $D_3$}

\textit{Dihedral groups} are groups of the symmetries of regular polygons. The first group we try is the simplest, $D_3$, the symmetries of an equilateral triangle. There are $6$ elements: $e, r, r^2, f, fr, fr^2$, where $f$ is a reflection, $r$ is a rotation, and $e$ is the identity. Here we map $e \mapsto 0$, $r \mapsto 1$, $r^2 \mapsto 2$, $f \mapsto 3$, $fr \mapsto 4$, $fr^2 \mapsto 5$.

We can generate a Cayley table of numbers to use for our conjugation relation:
\begin{equation*}
\begin{bmatrix}
0 & 1 & 2 & 3 & 4 & 5\\
1 & 0 & 3 & 2 & 5 & 4\\
2 & 4 & 0 & 5 & 1 & 3\\
3 & 5 & 1 & 4 & 0 & 2\\
4 & 2 & 5 & 0 & 3 & 1\\
5 & 3 & 4 & 1 & 2 & 0
\end{bmatrix}
\end{equation*}

\subsection{Coloring trefoil}
This is the diagram for the trefoil knot.
\begin{figure}[H]
  \centering
  \includegraphics[width=0.6\textwidth]{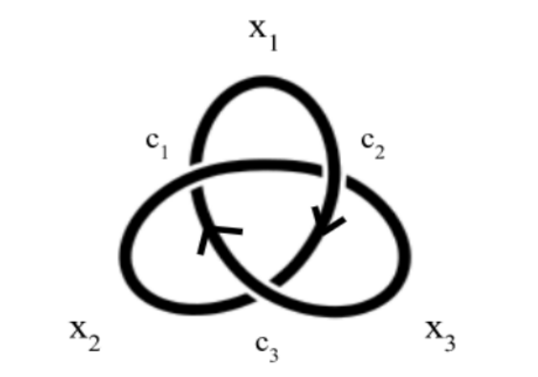}
  \caption{Crossing Relations [11]}
  \label{fig:diagram}
\end{figure} 

 We follow the arrows to label the crossings and the arcs for the diagram, giving us a table with the fundamental quandle coloring for the figure.
 
\begin{figure} [H]
    \centering
    \includegraphics[width=0.5\linewidth]{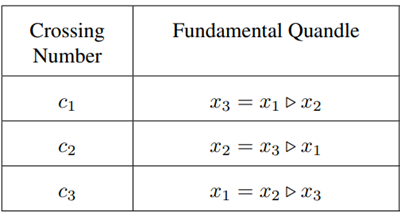}
    \caption{Fundamental Quandle Relations}
    \label{fig:placeholder}
\end{figure}

For the conjugation quandle over a dihedral group $D_3$, we obtain the relations,
\begin{align*}
\text{Crossing 1: }& x_3 = x_2^{-1} x_1 x_2 \quad (\text{over } D_3)\\
\text{Crossing 2: }& x_2 = x_1^{-1} x_3 x_1 \quad (\text{over } D_3)\\
\text{Crossing 3: }& x_1 = x_2^{-1} x_3 x_2 \quad (\text{over } D_3)
\end{align*}
which corresponds to this system of equations,
\begin{align*}
0 &= x_2^{-1} x_1 x_2 - x_3 \quad (\text{over } D_3)\\
0 &= x_1^{-1} x_3 x_1 - x_2 \quad (\text{over } D_3)\\
0 &= x_2^{-1} x_3 x_2 - x_1 \quad (\text{over } D_3)
\end{align*}

In order to solve this system of equations, we map each components to each element in the group, and take all the valid combinations to be the number of Homomorphisms.

https://www.wolframcloud.com/obj/zncubefan/Published/trefoilconjquandle.nb

We look for the counting invariant by assigning $x_1$, $x_2$, and $x_3$ to the elements in the dihedral group $D_3$, and see whether the equations hold for all assignments of $x$ to the group. Then we count the number of valid assignments. In this case, the valid assignments are \{0,0,0\}, \{1,1,1\}, \{1,2,5\}, \{1,5,2\}, \{2,1,5\}, \{2,2,2\}, \{2,5,1\}, \{3,3,3\}, \{4,4,4\}, \{5,2,1\}, \{5,1,2\}, \{5,5,5\}. There are $12$ valid homomorphisms in total, with $6$ monochromatic colorings, and $6$ solutions with $3$ elements. This means that the enhanced counting polynomial for the conjugation quandle of a trefoil over $D_3$ is $6q^3 + 6q$.

\subsection{Coloring Connected Sum of Hopf links}

We do the same process for the connected sum of the Hopf links.

\medskip
\begin{figure}[h]
    \centering
    \includegraphics[width=0.7\linewidth]{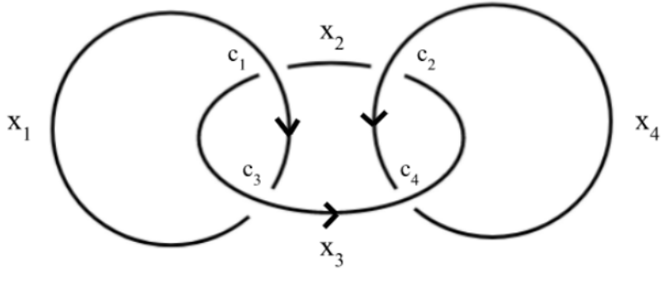}
    \caption{Colored Connected Sum of Hopf Links}
    \label{fig:placeholder}
\end{figure}

\medskip
For the conjugation quandle $y^{-1}xy$ corresponding to each fundamental quandle $x \triangleright y$ of $H\#H$ at each crossing over a dihedral group $D_3$, we obtain the relations
\begin{align*}
\text{Crossing 1: }& x_2 = x_1^{-1} x_3 x_1 \quad (\text{over } D_3)\\
\text{Crossing 2: }& x_3 = x_4^{-1} x_2 x_4 \quad (\text{over } D_3)\\
\text{Crossing 3: }& x_1 = x_3^{-1} x_1 x_3 \quad (\text{over } D_3)\\
\text{Crossing 4: }& x_4 = x_3^{-1} x_4 x_3 \quad (\text{over } D_3)
\end{align*}
which corresponds to this system of equations,
\begin{align*}
0 &= x_1^{-1} x_3 x_1 - x_2 \quad (\text{over } D_3)\\
0 &= x_4^{-1} x_2 x_4 - x_3 \quad (\text{over } D_3)\\
0 &= x_3^{-1} x_1 x_3 - x_1 \quad (\text{over } D_3)\\
0 &= x_3^{-1} x_4 x_3 - x_4 \quad (\text{over } D_3)
\end{align*}

Running the code in the same format as the trefoil, 

https://www.wolframcloud.com/obj/zncubefan/Published/hopfconjugationDn.nb  

We have $48$ total solutions, with $6$ containing $1$ coloring, $36$ containing $2$ colorings, and $6$ containing $3$ colorings, giving us a counting polynomial of $6q^3 + 36q^2 + 6q$.

This is the basis that we have to compare to the Allen--Swenberg link.

\subsection{Coloring the First Allen--Swenberg link}

Allen and Swenberg constructed a series of links that cannot be distinguished from the $H\#H$ via the Alexander-Conway polynomial. 

The simplest of which is their first link. 
\begin{figure}[H]
    \centering
    \includegraphics[width=0.5\linewidth]{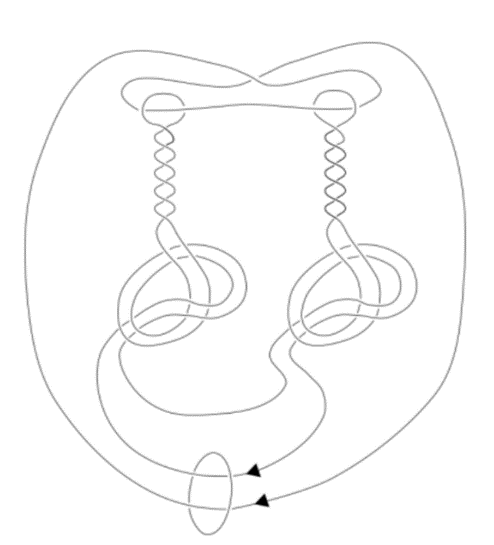}
    \caption{First Allen-Swenberg Link}
    \label{fig:placeholder}
\end{figure}
They showed that this link is different from $H\#H$ because their Jones polynomials are different.

Since the goal of this paper is to examine whether the conjugation quandle can detect causality along with the Alexander--Conway polynomial, our first task is to calculate the counting invariant, and the enhanced counting polynomial for this link to find if it can distinguish between the two.

First, label the diagram

\begin{figure} [H]
    \centering
    \includegraphics[width= 1\linewidth]{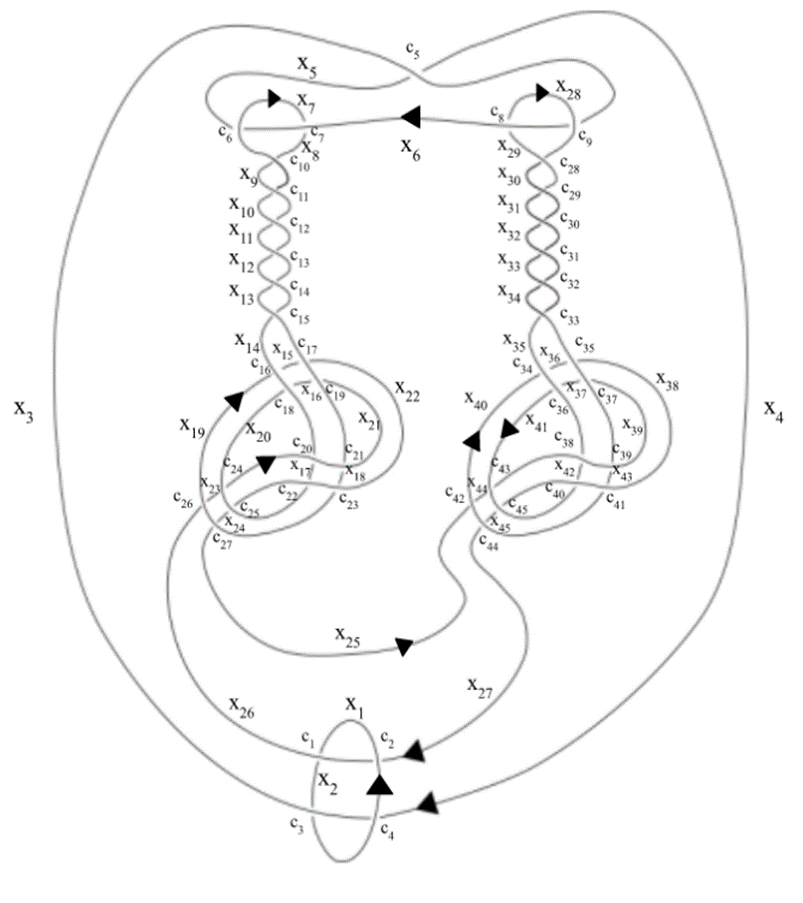}
    \caption{Colored Allen-Swenberg Link [11]}
    \label{fig:placeholder}
\end{figure}

Then, generate the Fundamental Quandle Table

\begin{figure} [H]
    \centering
    \includegraphics[width=0.8\linewidth]{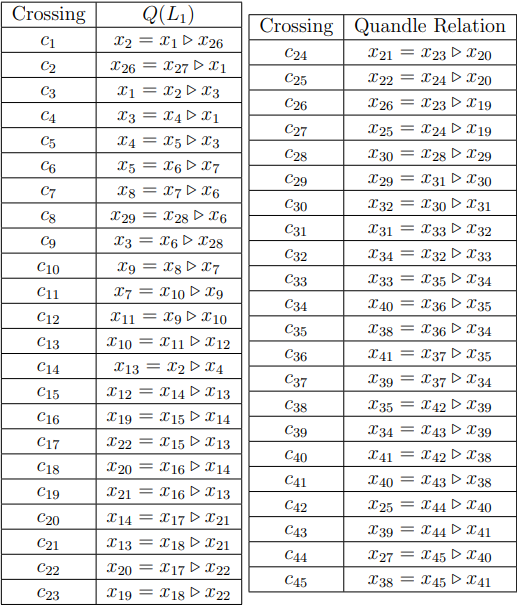}
    \caption{Fundamental Quandle Table for First Allen-Swenberg Link [9]}
    \label{fig:placeholder}
\end{figure}

Here, since everything is labeled, we can set the expression 
\[
x \triangleright y \;=\; y^{-1} x y
\]
and apply this for all of the crossings. Now we cycle across the dihedral group $D_3$ to find all of the valid solutions.

We use a different algorithm to run this link. If we tried to brute force the link, it would take $6^{45}$ combinations. Instead, the algorithm breaks the link down into sections, and runs each section independently, taking a divide and conquer approach. We then join the valid solutions together at the end.

https://www.wolframcloud.com/obj/zncubefan/Published/AllenSwenbergDn.nb 

After running a divide and conquer algorithm, we obtain a total of $48$ valid colorings, with $6$ being monochromatic, $36$ containing $2$ colors, and $6$ containing $3$ colors, giving us a counting polynomial of $6q^3 + 36q^2 + 6q$. This result is the same result as the result we obtained from $H\#H$ through the conjugation quandle in $D_3$. Therefore, the conjugation quandle in $D_3$ cannot show causality when coupled with the Alexander--Conway polynomial.

\subsection{Dihedral group $D_4$}

\sloppy The dihedral group $D_4$ is the permutation group for a square. The elements are $e, r, r^2, r^3, f, fr, fr^2, fr^3$. We map $e \mapsto 0$, $r \mapsto 1$, $r^2 \mapsto 2$, $r^3 \mapsto 3$, $f \mapsto 4$, $fr \mapsto 5$, $fr^2 \mapsto 6$, $fr^3 \mapsto 7$, giving a labeling of $D_4$ by $\{0,1,2,3,4,5,6,7\}$. We get the Cayley table
\begin{equation*}
\begin{bmatrix}
0&1&2&3&4&5&6&7\\
1&2&3&0&5&6&7&4\\
2&3&0&1&6&7&4&5\\
3&0&1&2&7&4&5&6\\
4&7&6&5&0&3&2&1\\
5&4&7&6&1&0&3&2\\
6&5&4&7&2&1&0&3\\
7&6&5&4&3&2&1&0
\end{bmatrix}
\end{equation*}
Now we can, again, calculate the conjugation quandle over $D_4$, for both $H\#H$, and the Allen--Swenberg link. This time, we obtain the result
\[72q^3 + 96q^2 + 8q\]
with $176$ total homomorphisms for both links. Therefore, we conclude that over the group of the permutations of a square, causality cannot be detected by using the conjugation quandle over $D_4$.

\subsection{Higher order groups}

The same calculation was done for $D_6$, and it was found to be $384$ homomorphisms, and a polynomial of $192q^3 + 180q^2 +12q$ for both links. For $D_7$, the polynomial was found to be $210q^3 + 168q^2 +14q$ with $392$ homomorphisms for both links. $D_8$ and $D_9$ were also computed, and while $D_9$ could distinguish the links, its size makes it computationally impractical, so we instead investigate $D_5$. 

One of the limitations of higher-order dihedral groups, such as $D_9$, is that they are computationally infeasible for practical purposes, making them less effective for detecting causality than $D_5$.

\subsection{Dihedral group $D_5$}

The dihedral group $D_5$ is the automorphism group of a regular pentagon. The elements are $e, r, r^2, r^3, r^4, f, fr, fr^2, fr^3, fr^4$. We map $e \mapsto 0$, $r \mapsto 1$, $r^2 \mapsto 2$, $r^3 \mapsto 3$, $r^4 \mapsto 4$, $f \mapsto 5$, $fr \mapsto 6$, $fr^2 \mapsto 7$, $fr^3 \mapsto 8$, $fr^4 \mapsto 9$. We get the Cayley table
\begin{equation*}
\begin{bmatrix}
0&1&2&3&4&5&6&7&8&9\\
1&2&3&4&0&9&5&6&7&8\\
2&3&4&0&1&8&9&5&6&7\\
3&4&0&1&2&7&8&9&5&6\\
4&0&1&2&3&6&7&8&9&5\\
5&6&7&8&9&0&1&2&3&4\\
6&7&8&9&5&4&0&1&2&3\\
7&8&9&5&6&3&4&0&1&2\\
8&9&5&6&7&2&3&4&0&1\\
9&5&6&7&8&1&2&3&4&0
\end{bmatrix}
\end{equation*}

Now we calculate the number of homomorphisms over the dihedral group $D_5$. For the $H\#H$, the enhanced quandle counting polynomial is
\[60q^3 + 90q^2 + 10q,\]
with $160$ total homomorphisms.

Calculating for the Allen--Swenberg link, we also find that there are a total of $200$ homomorphisms, and the counting polynomial
\[20q^8 + 20q^6 + 60q^3 + 90q^2 + 10q.\]

\textbf{Finally, we have a different result for $D_5$.}

For the group $D_5$, both the counting invariant and the counting polynomial are distinct for $H\#H$, and the Allen--Swenberg link. After trying all the smaller dihedral groups that are computable, only $D_5$ is able to detect causality.

\section{Results for the First Allen--Swenberg Link}

This table shows the results of the computations for the conjugation quandle over each dihedral group. 

\begin{figure} [H]
    \centering
    \includegraphics[width=1\linewidth]{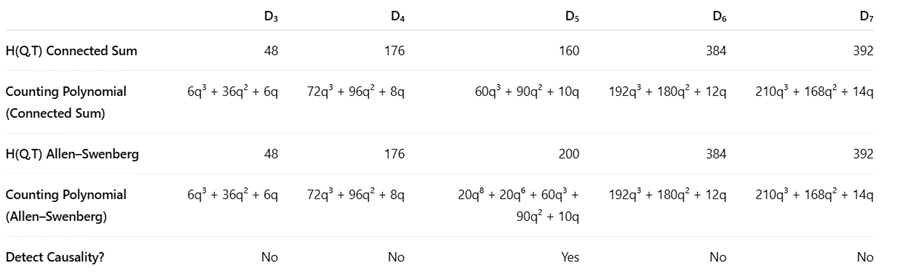}
    \caption{Results for Each Group}
    \label{fig:placeholder}
\end{figure}

The conjugation quandle across the group $D_5$ can detect causality in links, when combined with the Alexander--Conway polynomial. However, when run across other dihedral groups, like $D_3$, $D_4$, $D_6$, and $D_7$ the conjugation quandle fails to detect causality.

What is interesting about this result, is that while running across other dihedral groups, even the enhanced counting polynomial failed to detect causality. However, for the dihedral group $D_5$, not only was the enhanced counting polynomial distinct, even the counting invariant was different.

\begin{theorem}\label{thm:main}
Let $L_1$ be the first Allen--Swenberg $2$-sky link and $H\#H$ be the connected sum of the $2$ Hopf links. Then $|\mathrm{Hom}(Q(L_1), \mathrm{Conj}(D_5))| = 200$ and $|\mathrm{Hom}(Q(H\#H), \mathrm{Conj}(D_5))| = 160$, and the enhanced counting polynomials differ. In particular, $L_1$ and $H\#H$ are distinguished by the $D_5$ conjugation quandle.
\end{theorem}

However, this result cannot prove that the conjugation quandle across $D_5$ can distinguish between infinite series of the Allen--Swenberg, therefore we must turn to the second Allen--Swenberg link.

\section{Generalized Results Across $D_5$}

The second Allen--Swenberg link is the following link.

\begin{figure} [H]
    \centering
    \includegraphics[width=0.9\linewidth]{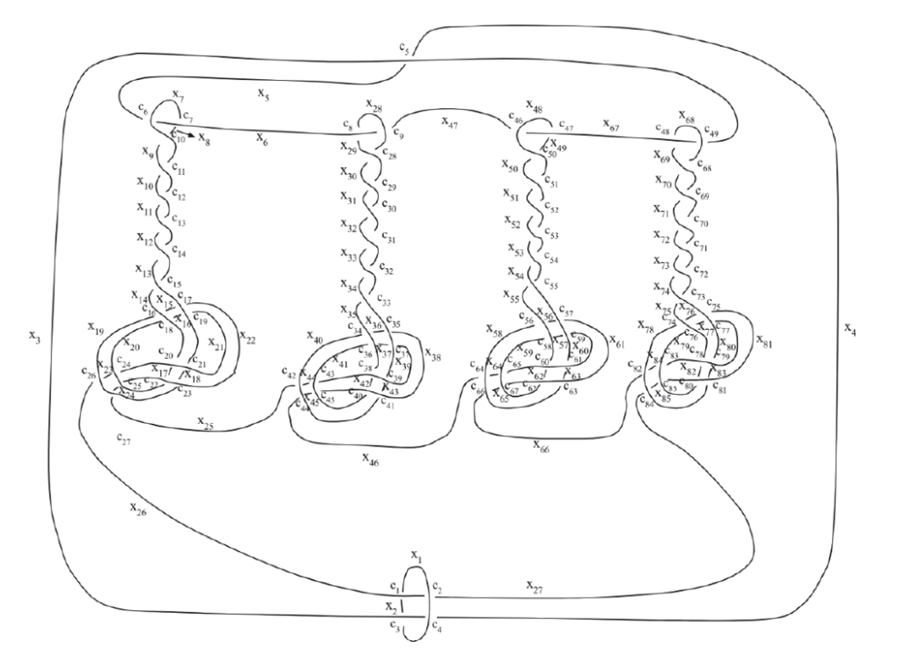}
    \caption{Colored Second Allen-Swenberg Link [9]}
    \label{fig:placeholder}
\end{figure}

Make a table for the crossings, where each crossing is related to the arcs coming out of the crossing:
\begin{align*}
\{26, 1, 2\}, \{1, 27, 26\}, \{3, 2, 1\}, ...
\end{align*}
where each triple represents the relation $\text{arc}_3 \triangleright \text{arc}_2 = \text{arc}_1$ at each crossing.After running the divide and conquer algorithm, 

https://www.wolframcloud.com/obj/zncubefan/Published/AllenSwenberg2Dn.nb

we obtain the following results. We find that the quandle counting invariant for $D_3$ distinguishes between $L_1$ and $L_2$, with $48$ homomorphisms for $L_1$, and $54$ for $L_2$, while $D_4$ cannot distinguish between $H\#H$ or $L_1$ and $L_2$.

Across the dihedral group $D_5$, we have seen that both the first and second Allen--Swenberg links $L_1$ and $L_2$ can be distinguished from the connected sums of the Hopf links $H\#H$ through the counting invariant alone. We have also found that in $D_5$, the counting polynomial can distinguish between $L_1$, and $L_2$, with $L_1$ having a counting polynomial of $20q^8 + 20q^6 + 60q^3 + 90q^2 + 10q$, while $L_2$ has a counting polynomial of $20q^7 + 20q^5 + 60q^3 + 90q^2 + 10q$.

We see that the higher order terms in the polynomial change from $8$ and $6$ to $7$ and $5$. This means that the number of colors used for the Allen--Swenberg links differ while the number of colorings remain the same.

We want to prove that any Allen--Swenberg link can be distinguished from the connected sum of the Hopf links.

\begin{figure} [H]
    \centering
    \includegraphics[width=0.8\linewidth]{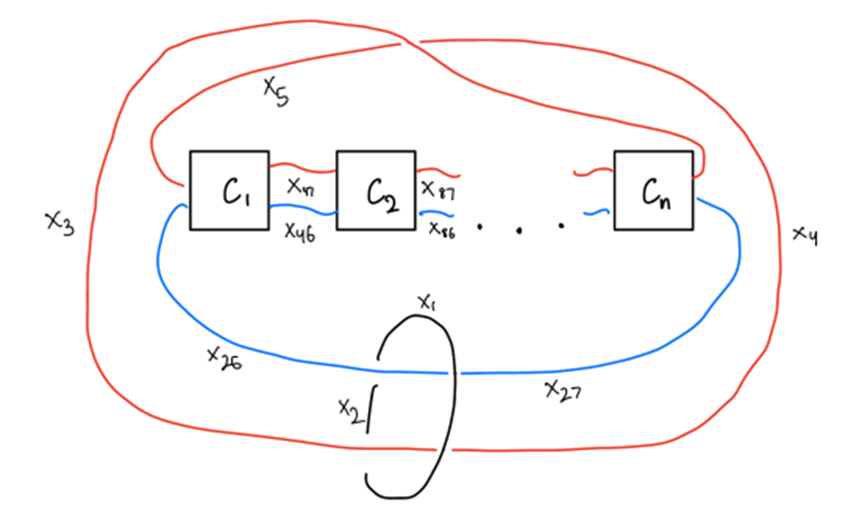}
    \caption{Generating New Allen-Swenberg Links [9]}
    \label{fig:placeholder}
\end{figure}

In order to generate a new Allen--Swenberg link, you take $C_1$, the section of the knot with crossings $10$--$26$, and $29$--$45$, and repeat the crossings to form a new section, $C_2$, with crossings $50$--$66$, and $69$--$85$. In order to generate $C_n$, you do the same process, all while setting $x_5 = x_{47} = x_{87} = \cdots = x_3$, and $x_{26} = c_{46} = c_{86} = \cdots = 27$.

\begin{lemma}\label{lem:parity} \sloppy
For every coloring $f_1 \in \mathrm{Hom}(Q(L_1), \mathrm{Conj}(D_5))$ there exists $f_2 \in \mathrm{Hom}(Q(L_2), \mathrm{Conj}(D_5))$, where there exists either $f_1 = f_2$, or a parity case where the number of homomorphisms stays the same, but the number of colors found in each homomorphism differs.
\end{lemma}

\begin{figure}[H]
    \centering
    \includegraphics[width=0.9\linewidth]{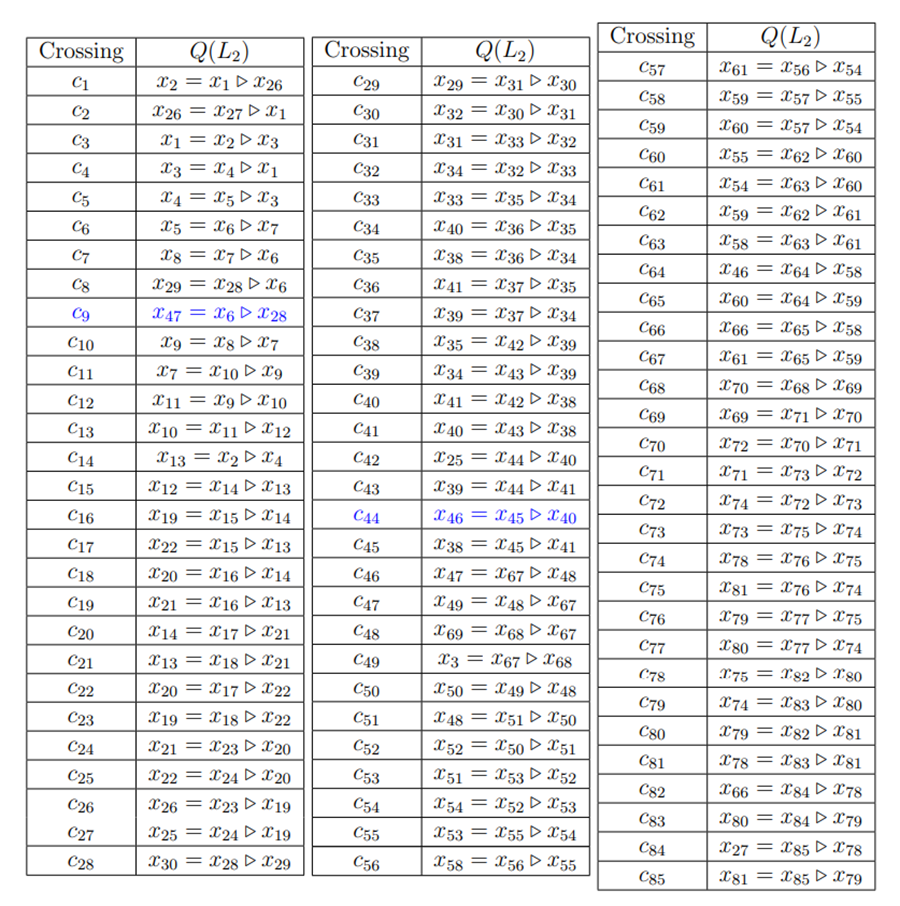}
    \caption{Table for Allen-Swenberg Second Link Colorings [9]}
    \label{fig:placeholder}
\end{figure}

From the table above, we find that the crossings $1$--$45$ found in the first link are copied in the second link, except for that of $C_9$ and $C_{44}$. Jain, in his paper, proved the first part of Lemma~\ref{lem:parity}. Our goal is to prove the second part.

\begin{figure} [H]
    \centering
    \includegraphics[width=0.9\linewidth]{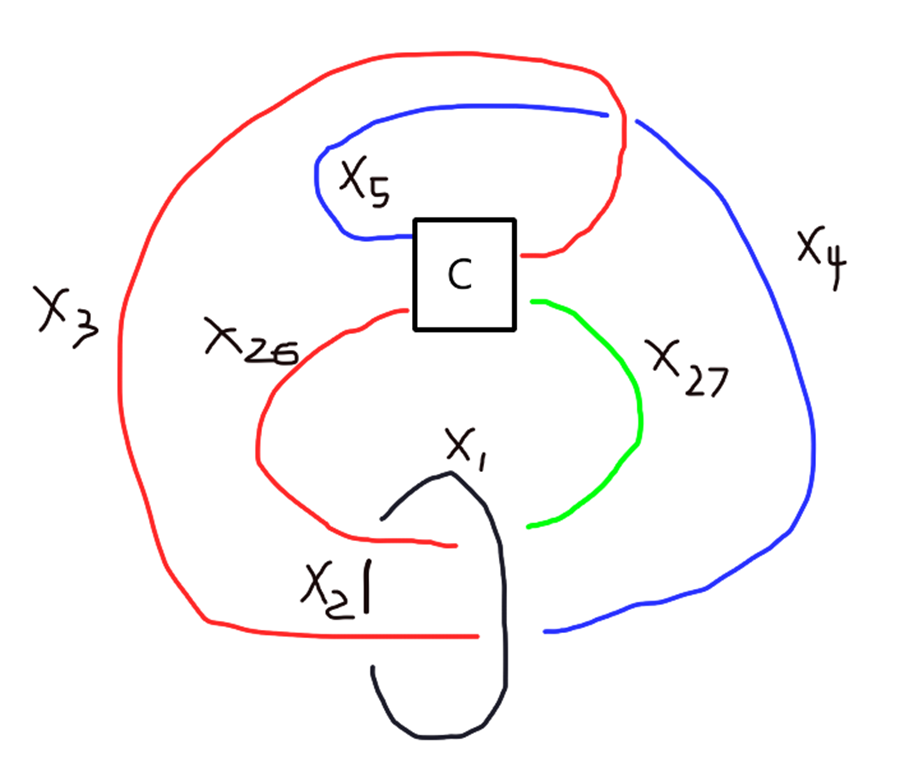}
    \caption{Hypothetical Coloring for $L_1$}
    \label{fig:placeholder}
\end{figure}

If we draw $L_1$ as such, we see that it contains one section, which we can name $C$. This is a hypothetical coloring for the first Allen--Swenberg link. Here we assume the black component of the link is an extra coloring when $x_{26}$ is red and when $x_{27}$ is green yields a new color for $x_1$ and $x_2$; black is a valid coloring. And hence we have $4$ valid colors outside $C$.

\begin{figure} [H]
    \centering
    \includegraphics[width=1\linewidth]{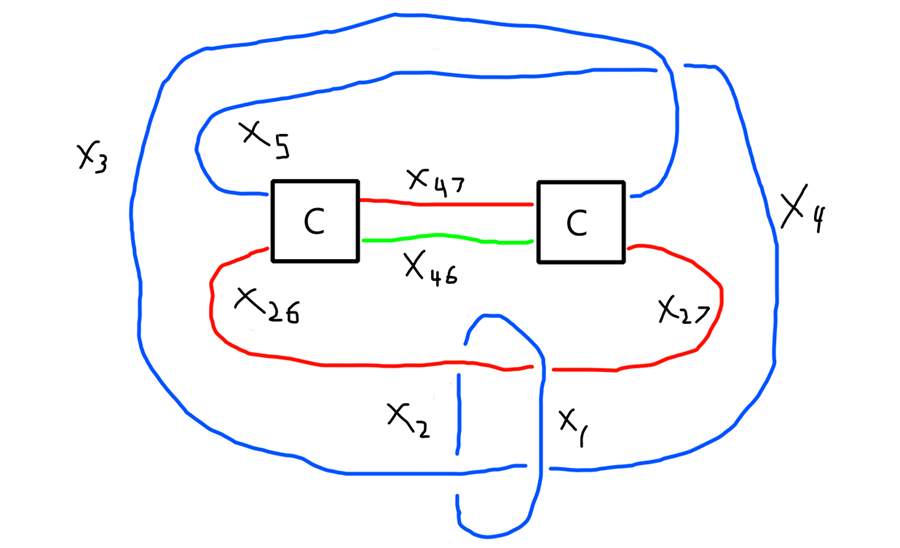}
    \caption{Hypothetical Coloring for $L_2$}
    \label{fig:placeholder}
\end{figure}

However, we see the parity case in the second Allen--Swenberg Link. Here, $[x_5, x_{26}]$, going into the crossing junction $C$, where the parity flips it so that $[x_{56}, x_{47}]$ is green, but it is flipped again, so that $[x_3, x_{27}]$ coming out of the junction is then blue and red, leading to the third component, which is comprised of $x_1$ and $x_2$, to be a different valid coloring. This link only has $4$ colorings outside of $C$.

This explains the change from $20q^8 + 20q^6 + 60q^3 + 90q^2 + 10q$ to $20q^7 + 20q^5 + 60q^3 + 90q^2 + 10q$. Here, the $6$ and $8$ colorings in $L_1$ run into the parity case, and therefore, are reduced by one coloring. This results in the new link having $5$ and $7$ colorings. The rest of the link does not run into the parity case, and therefore, the $3$, $2$, and trivial coloring polynomials do not decrease in size.

This means that the first and odd links of the infinite series will have the same parity, while the even ones will have the same parity.

Since Jain proved that every coloring for each parity also exists in the next link, we prove that
$\left|\mathrm{Hom}(\text{Allen--Swenberg Series}, \mathrm{Conj}(D_5))\right| \ge 200$, 
which is greater than 
$|\mathrm{Hom}(H\#H, \mathrm{Conj}(D_5))|$.

Therefore,

\begin{theorem}\label{thm:series}
The conjugation quandle across $D_5$ can distinguish between all Allen--Swenberg $2$-sky links, and the connected sum of the Hopf links.
\end{theorem}

\begin{conjecture}\label{conj:dihedral-collection}
We conjecture the collection of all conjugation quandles over all dihedral groups is able to detect causality across $(2+1)$ globally-hyperbolic spacetimes when coupled with the Alexander--Conway polynomial.
\end{conjecture}

\section{Conclusion}
We have clearly shown that the first and second Allen--Swenberg link can be differentiated from the connected sum of the Hopf links by the conjugation quandle over $D_5$. We have also proved that the conjugation quandle over $D_5$ can differentiate between all Allen--Swenberg links, and therefore, when coupled with the Alexander--Conway polynomial, can determine causality between events in a $(2+1)$-dimensional globally hyperbolic spacetime. However, we have not tested other non-abelian groups, which are left open for further research.

\section*{Acknowledgements}
This project was completed as part of the Summer 2025 Theoretical Mathematics and Knot Theory Seminar by the Horizon Academic Research Program. It was completed under the supervision and guidance of Professor Vladimir Chernov from Dartmouth College, and Dr.\ Ryan Maguire from MIT. I would like to thank Professor Chernov and Dr.\ Maguire for their commitment to this class.


\begin{thebibliography}{99}
\bibitem{AllenSwenberg2021}
S.~Allen and J.~Swenberg, \emph{Do link polynomials detect causality in globally hyperbolic spacetimes?} J.~Math.~Phys., 62(3), 2021.

\bibitem{BernalSanchez2003}
A.~Bernal and M.~S\'anchez, \emph{On smooth Cauchy hypersurfaces and Geroch's splitting theorem}, Commun.~Math.~Phys., 243:461--470, 2003.

\bibitem{BernalSanchez2007}
A.~Bernal and M.~S\'anchez, \emph{Globally hyperbolic spacetimes can be defined as ``causal'' instead of ``strongly causal''}, Class.~Quant.~Grav., 24:745--750, 2007.

\bibitem{Chen2024}
H.~Chen, \emph{Medial quandles' capability of detecting causality and properties of their coloring on certain links and knots}. arXiv:2411.04477v4 [math.GT], 27 Nov 2024. \url{https://arxiv.org/pdf/2411.04477}

\bibitem{ChernovNemirovski2016}
V.~Chernov and S.~Nemirovski, \emph{Causality and Legendrian linking for higher dimensional spacetimes}. Journal of Geometry and Physics 133, 26--29. J.~Symplectic Geom.~14 (2016) 149--170.

\bibitem[5.1]{ChernovNemirovski-Low} V. Chernov and S. Nemirovski, \emph{Legendrian links, causality, and the Low conjecture}, arXiv:0810.5091 [math.SG] (submitted Oct 2008); Geom. Funct. Anal. 19 (2010), no. 5, 1320–1333.

\bibitem[5.2]{ChernovNemirovski-NonNeg} V. Chernov and S. Nemirovski, \emph{Non-negative Legendrian isotopy in $ST^*M$}, arXiv:0905.0983 [math.SG] (submitted May 2009); Geom. Topol. 14 (2010), no. 1, 611–626.

\bibitem{ChernovMartinPetkova2020}
V.~Chernov, G.~Martin, and I.~Petkova, \emph{Khovanov homology and causality in spacetimes}. J.~Math.~Phys., 61(0222503), 2020.

\bibitem{Geroch1970}
R.~Geroch, \emph{Domain of dependence}. J.~Math.~Phys., 11:437--449, 1970.

\bibitem{HawkingEllis1973}
S.~W.~Hawking and G.~F.~R.~Ellis, \emph{The large scale structure of space-time}. Cambridge Monographs on Mathematical Physics 1. Cambridge University Press (1973).

\bibitem{Jain2024}
A.~Jain, \emph{Detecting causality with symplectic quandles}, Lett.\ Math.\ Phys., 114:3 (2024). arXiv:2310.06853v1 [math.GT].

\bibitem{Joyce1982}
D.~Joyce, \emph{A classifying invariant of knots, the knot quandle}. J.~Pure Appl.~Algebra, 23:37--65, 1982.

\bibitem{Leventhal2023}
J.~Leventhal, \emph{Alexander quandles and detecting causality}. arXiv, 2209.05670v1, 2023.

\bibitem{Low1988}
R.~J.~Low, \emph{Causal relations and spaces of null geodesics}. PhD thesis, Oxford University, 1988.

\bibitem{Matveev1984}
S.~V.~Matveev, \emph{Distributive groupoids in knot theory}. Math. USSR-S, 47(7383), 1984.

\bibitem{NatarioTod2004}
J.~Nat\'ario and P.~Tod, \emph{Linking, Legendrian linking and causality}. Proc.\ London Math.\ Soc., 88:251--272, 2004.

\bibitem{NelsonNotes2004}
S.~Nelson, \emph{Quandles and racks}. 2004. \url{https://www1.cmc.edu/pages/faculty/WNelson/quandles.html}.

\bibitem{Penrose1979}
R.~Penrose, \emph{Singularities and time-asymmetry}. 1979. In \emph{General Relativity: An Einstein Centenary Survey}, eds. S.~W.~Hawking \& W.~Israel, pp.~581--638. Cambridge University Press.

\bibitem{Wolfram2025}
I.~Wolfram Research. \emph{Mathematica desktop}, version 14.3, 2025.
\end{thebibliography}
\end{document}